\documentclass{amsart}
\usepackage{graphicx}
\usepackage{amssymb}
\usepackage{amsmath}
\newtheorem{theorem}{Theorem}[section]
\newtheorem{lemma}{Lemma}[section]

\newenvironment{pol}{{\bf Proof of Lemma 3.1.}}{\hfill\fbox{}\par\vspace{.2cm}}
\newenvironment{pot1}{{\bf Proof of Theorem 1.5.}}{\hfill\fbox{}\par\vspace{.2cm}}
\newenvironment{pot2}{{\bf Proof of Theorem 1.6.}}{\hfill\fbox{}\par\vspace{.2cm}}
\numberwithin{equation}{section}

\def\charf {\mbox{{\text 1}\kern-.24em {\text l}}}

\def\bea{\begin{eqnarray*}}
\def\eea{\end{eqnarray*}}
\def\be{\begin{eqnarray}}
\def\ee{\end{eqnarray}}

\begin{document}

\title{Integral Curvature Bounds and Bounded Diameter with Bakry--Emery Ricci Tensor}


\author{Seungsu Hwang}
\address{Department of Mathematics, Chung-Ang University, 84 Heukseok-ro, Dongjak-gu, Seoul, Republic of Korea}
\curraddr{}
\email{seungsu@cau.ac.kr}
\thanks{}

\author{Sanghun Lee}
\address{Department of Mathematics, Chung-Ang University, 84 Heukseok-ro, Dongjak-gu, Seoul, Republic of Korea}
\curraddr{}
\email{kazauye@cau.ac.kr}
\thanks{}


\date{}

\dedicatory{}

\begin{abstract}
 For Riemannian manifolds with a smooth measure $(M, g, e^{-f}dv_{g})$, we prove a generalized Myers compactness theorem when Bakry--Emery Ricci tensor is bounded from below and $f$ is bounded.
\end{abstract}

\maketitle
\section{Introduction}
One of the most fundamental results in Riemannian geometry is the Myers theorem, which states that if a complete Riemannian manifold $(M,g)$ satisfies $Ric \geq (n-1)H $ with $H>0$, then $M$ is compact and ${\rm diam}(M) \leq \frac{\pi}{\sqrt{H}}$. Here, $Ric$ is the Ricci curvature of the metric $g$. This theorem has been generalized through different approaches (see \cite{AB}, \cite{GL}, \cite{CS}, and \cite{WR}), one of which is the effort of Wei and Wylie, who proved the theorem for manifolds with a positive lower Bakry-–Emery Ricci curvature bound in \cite{GW}. A Bakry-–Emery Ricci tensor is defined as
$$
Ric_{f} = Ric+ {\rm Hess}\, f,
$$
where $f$ is a smooth function on $M$ and ${\rm Hess}\, f$ is the hessian of $f$.
Where $|f|\leq k$, they proved that
$$
{\rm diam}(M) \leq \frac{\pi}{\sqrt{H}} + \frac{4k}{(n-1)\sqrt{H}}.
$$
Several works have attempted to generalize this result (for example, see \cite{MR}, \cite{YS}, \cite{JY}, and \cite{SZ}), including that of Sprouse which can be summarized in the following three theorems.
\begin{theorem} [\cite{CS}]
 Let $(M,g)$ be a compact Riemannian manifold of nonnegative Ricci curvature. Then, for any $\delta > 0$, there exists $\epsilon = \epsilon(n,\delta)$ such that if
\[ \frac{1}{vol(M)} \int_{M} ((n-1) - Ric_{-})_{+} dv_g < \epsilon(n,\delta), \]
then ${\rm diam}(M) < \pi + \delta$.
\end{theorem}
Here, $dv_{g}$ is the Riemannian volume density on $M$, $Ric_{-}(x)$ is the lowest eigenvalue of the Ricci tensor $Ric(x)$, and $h_+(x)=\max\{h(x), 0\}$ for an arbitrary function $h$ on $M$. For $Ric \geq (n-1)k$ with $k\leq 0$, these generalizations are attained.  
\begin{theorem}  
[\cite{CS}]
 Let $(M,g)$ be a complete Riemannian manifold with $Ric \geq (n-1)k$, $k \leq 0$. Then, for any $R, \delta > 0$, there exists $\epsilon = \epsilon(n, k, R, \delta)$ such that if
\[ \sup_{x} \frac{1}{vol(B(x,R))} \int_{B(x,R)} ((n-1) - Ric_{-})_{+} dv_g < \epsilon(n, k, R, \delta), \]
then $(M,g)$ is compact and ${\rm diam}(M) < \pi + \delta$.
\end{theorem}
 
\begin{theorem} [\cite{CS}]
 Let $(M,g)$ be a complete Riemannian manifold with $Ric \geq (n-1)k$, $k \leq 0$. Then, for any $R > 0$, there exists $\tilde{\epsilon} = \tilde{\epsilon}(n, k, R)$ such that if
\[ \sup_{x} \frac{1}{vol(B(x,R))} \int_{B(x,R)} ((n-1) - Ric_{-})_{+} dv_g < \tilde{\epsilon}(n, k, R), \]
then the universal cover of $M$ is compact, and hence, $\pi_{1}(M)$ is finite.
\end{theorem}

 We will generalize these results to 
the Bakry-–Emery Ricci tensor bounded from below. Let $(M,g,e^{-f}dv_g)$ be a smooth metric measure space, where $M$ is a complete $n$-dimensional Riemannian manifold with metric $g$. Likewise, let $Ric_{f-}(x)$ denote the lowest eigenvalue of the Bakry-–Emery Ricci tensor $Ric_{f}(x)$. Then, we prove the following theorem.

\begin{theorem}\label{main1}
 Let $(M,g,e^{-f}dv_g)$ be a compact $n$-dimensional Riemannian manifold with $Ric_f \geq 0$ and $|f| \leq k$. Then, for any $\delta > 0$, there exists $\epsilon = \epsilon(n+4k,\delta)$ such that if
$$
\frac{1}{vol_f(M)} \int_{M} ((n-1) - Ric_{f-})_+ e^{-f}dv_g < \epsilon(n+4k,\delta),
$$
then ${\rm diam}(M) < \pi + \delta$.
\end{theorem}
   Given that $(M,g)$ is noncompact or does not exhibit a nonnegative Bakry–-Emery Ricci curvature, averaging the bad part of $Ric_f$ over metric ball, as in \cite{CS} yields a similar result as follows.
\begin{theorem}\label{main2}
 Let $(M,g,e^{-f}dv_g)$ be a complete $n$-dimensional Riemannian manifold with $Ric_{f} \geq (n-1)H$, $H < 0$, and $|f| \leq k$. Then, for any $R, \delta > 0$, there exists $\epsilon = \epsilon(n+4k, H, R, \delta)$ such that if
$$
\sup_{x} \frac{1}{vol_{f}(B(x,R))} \int_{B(x,R)} ((n-1) - Ric_{f-})_{+} e^{-f}dv_{g} < \epsilon(n+4k,H,R,\delta),
$$
then $M$ is compact and ${\rm diam}(M) < \pi + \delta$.
\end{theorem}
Finally, we could obtain the result for the fundamental group of $M$, which is stated as follows.

\begin{theorem}\label{main3}
 Let $(M,g,e^{-f}dv_{g})$ be a complete $n$-dimensional Riemannian manifold with $Ric_{f} \geq (n-1)H$, $H<0$, and $|f| \leq k$. Then for any $R>0$, there exists $\tilde{\epsilon} = \tilde{\epsilon}(n+4k,H,R)$ such that if
$$
\sup_{x} \frac{1}{vol_{f}(B(x,R))} \int_{B(x,R)} ((n-1) - Ric_{f-})_{+} e^{-f}dv_{g} < \tilde{\epsilon}(n+4k,H,R),
$$
then the universal cover of $M$ is compact, and hence, $\pi_{1}(M)$ is finite.
\end{theorem}

\section{Proof of Theorem~\ref{main1} }
Let $(M,g,e^{-f}dv_g)$ be a smooth metric measure space, where $(M,g)$ is a complete $n$-dimensional Riemannian manifold.
 Let $A_1, A_2, W$ be open subsets of $M$ such that $A_1,A_2 \subset W$, and all minimal geodesics $\gamma_{x,y}$ from $x \in A_1$ to $y \in A_2$ lie in $W$.
 
 We will use the estimate of Cheeger and Colding for Bakry–Emery Ricci tensor (\cite{ML}, Proposition 2.3); thus, for a nonnegative integrable function on $M$, 
\begin{eqnarray}
\lefteqn{\int_{A_1 \times A_2} \int_{\gamma_{x,y}} h(\gamma(s)) \, ds(e^{-f}dv_g)^{2}}\nonumber\\
& & \leq  C(n+4k,H,R) ({\rm diam}(A_{2})vol_{f}(A_1) + {\rm diam}(A_{1})vol_{f}(A_2))\nonumber\\
& & \times \int_{W} h e^{-f}dv_g, \label{eqn1}\end{eqnarray}
then,
$$
C(n+4k,H,R) = \sup_{0 < \frac{s}{2} \leq u \leq s} \frac{A_{H}^{n+4k}(s)}{A_{H}^{n+4k}(u)},$$
and$$
R \geq \sup \{d(x,y) \mid (x,y) \in (A_{1} \times A_{2})\},$$
where $A_{H}^{n+4k}(r)$ denotes the area element on $\partial B(r)$ in $M_H^{n+4k}$, the simply connected model space of dimension $n+4k$ with constant curvature $H$.
Because $H=0$, we denote $C(n+4k,H,R)$ by  $C(n+4k)$. 
 
 Applying $sn_{H}(r)$ as a solution to $$sn_{H}'' + H sn_{H} = 0$$, $sn_{H}(0) = 0$ and $sn_{H}'(0) = 1$ are satisfied. Moreover, if $H=0$, then $$m_{H}^{n} = (n-1)\frac{sn_{H}'}{sn_{H}},$$ with the solution $sn_{H}(r) = r$.
 
   By the mean curvature comparison (3.15) in the proof of Theorem 1.1 in \cite{GW}, we have
\begin{eqnarray}
sn_{H}^{2}(r)m_{f}(r) &\leq & sn_{H}^{2}(r)m_{H}(r) - f(r)(sn_{H}^{2}(r))' \nonumber\\
& & + \int_{0}^{r} f(t)(sn_{H}^{2})''(t) dt. \label{eqn2}
\end{eqnarray}
Thus,
\[m_{f}(r) \leq m_{H}(r) + \frac{4k}{sn_{H}(r)}, \]
implying that
$$
m_{f}(r) \leq m_{H}(r)\left(1 + \frac{4k}{n-1}\right) = m_{H}^{n+4k}(r).
$$
Therefore,
\be
\frac{vol_{f}(B(p,R))}{vol_{f}(B(p,r))} \leq \frac{vol_{H}^{n+4k}(R)}{vol_{H}^{n+4k}(r)}.
\label{eqn3}\ee
 \vskip .5pc
   Now, let $p,q \in M$ such that $d(p,q) = {\rm diam}(M) = D$, $r>0$, $A_{1} = B(p,r)$, and $A_{2} = B(q,r)$. Applying the inequality (\ref{eqn1}),
\bea
\lefteqn{\int_{A_1 \times A_2} \int_{\gamma_{x,y}} ((n-1) - Ric_{f-})_{+} ds(e^{-f}dv_g)^{2} }\\
&  \leq &C(n+4k) (2r \, vol_{f}(A_1) + 2r \, vol_{f}(A_2))\int_{M} ((n-1) - Ric_{f-})_{+} e^{-f}dv_g.
\eea
   Consequently, let $vol_{H}^{n+4k}(r)$ be the volume of the radius $r$-ball  in $M_{H}^{n+4k}$, the simply connected model space of dimension $n+4k$ with constant curvature $H$. Then, we have
\bea
\lefteqn{\inf_{(x,y) \in A_{1} \times A_{2}} \int_{\gamma_{x,y}} ((n-1) - Ric_{f-})_{+} ds}\\
&\leq &2r \, C(n+4k) (\frac{1}{vol_{f}(A_{1})} + \frac{1}{vol_{f}(A_{2})}) \int_{M}((n-1) - Ric_{f-})_{+} e^{-f}dv_g\\ 
&\leq & 4r \, C(n+4k) \frac{vol_{H}^{n+4k}(D)}{vol_{H}^{n+4k}(r)}\frac{1}{vol_f(M)} \int_{M}((n-1) - Ric_{f-}) e^{-f}dv_g,\eea
where the last inequality follows from (\ref{eqn3}).

   Note that if $H = 0$, the volume element $v(r) = r^{n+4k-1}$, which gives
$$
vol_{H}^{n+4k}(r) = \int_{S^{n+4k-1}} ds^{n+4k-1} \int_{0}^{r} t^{n+4k-1} dt,
$$
to obtain
$$
\frac{vol_{H}^{n+4k}(D)}{vol_{H}^{n+4k}(r)} = \frac{D^{n+4k}}{r^{n+4k}}.
$$
Therefore,
\begin{eqnarray}
\lefteqn{\inf_{(x,y) \in (A_{1} \times A_{2})} \int_{\gamma_{x,y}} ((n-1) - Ric_{f-})_{+} ds}\nonumber \\
& \leq& 4r \, C(n+4k) \frac{D^{n+4k}}{r^{n+4k}} \frac{1}{vol_{f}(M)} \int_{M} ((n-1) -Ric_{f-})_{+} e^{-f}dv_g. \label{eqn4}
\end{eqnarray}

   Now, we can find a minimizing unit speed geodesic $\gamma$ from $x \in \overline{A_{1}}$ to $y \in \overline{A_{2}}$ of length $L = d(x,y)$.
Let $\{E_{1},\cdots, E_{n} = \gamma'\}$ be a parallel orthonormal frame along $\gamma$ and a smooth function $b \in C^{\infty}([0,L])$ such that $b(0) = b(L) = 0$; then, by the second variation of $\gamma$, we have
\bea
\sum_{i=1}^{n-1} I(bE_{i},bE_{i}) &=& \int_{0}^{L} (b')^{2}(n-1)\, dt - \int_{0}^{L} b^{2} Ric_{f}(\gamma', \gamma') \, dt \\
& & + \int_{0}^{L}b^{2}Hess(f)(\gamma',\gamma') \, dt.
\eea
Likewise, note that \[\int_{0}^{L} b^{2}Hess(f)(\gamma',\gamma') dt = \int_{0}^{L} b^{2} \frac{d}{dt} \langle \nabla f, \gamma'\rangle \, dt,\]
thus,
\bea
\int_{0}^{L} b^{2} \frac{d}{dt} \langle \nabla f,\gamma' \rangle dt &=& \int_{0}^{L}\left( -2bb' \frac{d}{dt}(f(\gamma(t))) + \frac{d}{dt}(b^{2}\langle \nabla f,\gamma'\rangle)\right) dt \\
&= &  \int_{0}^{L} \left( 2f\frac{d}{dt}(bb') - 2\frac{d}{dt}(fbb') + \frac{d}{dt}(b^{2}\langle \nabla f,\gamma'\rangle)\right) dt \leq 0,
\eea
such that,
$$
\sum_{i=1}^{n-1} I(bE_{i},bE_{i}) \leq \int_{0}^{L} (b')^{2}(n-1) dt - \int_{0}^{L} b^{2} Ric_{f}(\gamma', \gamma')dt.
$$
If we set the function $b$ as $b(t)=\sin\left(\frac{\pi t}{L}\right)$, then we obtain $(b'(t))^{2} = \frac{\pi^{2}}{L^{2}} \cos^{2}(\frac{\pi t}{L})$ and $b^{2}(t) = \sin^{2}\left(\frac{\pi t}{L}\right)$. Thus, 
\bea
\sum_{i=1}^{n-1} I(bE_{i},bE_{i}) &\leq & (n-1)\int_{0}^{L} \frac{\pi^{2}}{L^{2}} \cos^{2}\left(\frac{\pi t}{L}\right)dt - \int_{0}^{L} \sin^{2}\left(\frac{\pi t}{L}\right)Ric_{f}(\gamma',\gamma')\, dt,\\
& =&  (n-1) \frac{\pi^{2}}{L^{2}}\int_{0}^{L}\cos^{2}\left(\frac{\pi t}{L}\right)dt - (n-1) \int_{0}^{L} \sin^{2}\left(\frac{\pi t}{L}\right)\, dt\\
& & + \int_{0}^{L} \sin^{2}\left(\frac{\pi t}{L}\right)((n-1) - Ric_{f}(\gamma',\gamma'))\, dt,\\
&=& -\frac{(n-1)L}{2} \left(1 - \frac{\pi^{2}}{L^{2}}\right) + \int_{0}^{L} \sin^{2}\left(\frac{\pi t}{L}\right)((n-1) - Ric_{f}(\gamma',\gamma'))dt\\
&\leq & -\frac{(n-1)L}{2} \left(1 - \frac{\pi^{2}}{L^{2}}\right) + \int_{0}^{L}((n-1) - Ric_{f-})_{+} dt. \eea
By the inequality (\ref{eqn4}),
\bea
\sum_{i=1}^{n-1} I(bE_{i},bE_{i}) &\leq & -\frac{(n-1)L}{2}\left(1 - \frac{\pi^{2}}{L^{2}}\right) \\
& & + 4rC(n+4k)\frac{D^{n+4k}}{r^{n+4k}} \frac{1}{vol_{f}(M)}\int_{M}((n-1) - Ric_{f-})_{+} e^{-f}dv_g.
\eea
Now, let $r = \frac{D}{N}$, and choose $N=N(\delta)$ such that 
\begin{equation}
\frac{1}{1 - \frac{2}{N}} < \frac{\pi + \delta}{\pi + \frac{\delta}{2}}, \label{eqn5}
\end{equation}
by the triangle inequality, 
\begin{equation}
L = d(x,y) \geq d(p,q) - d(p,x) - d(y,q) = D\left(1 - \frac{2}{N}\right).\label{eqn6}
\end{equation}
Thus,
\bea
\sum_{i=1}^{n-1} I(bE_{i},bE_{i}) &\leq &-\frac{(n-1)L}{2}(1 - \frac{\pi^{2}}{L^{2}}) \\
& & + 4C(n+4k)\frac{L}{1 - \frac{2}{N}} \frac{N^{n+4k-1}}{vol_{f}(M)}\int_{M}((n-1) - Ric_{f-})_{+}e^{-f}dv_{g}.\eea
Setting
$$
\epsilon = \frac{(n-1)(1-\frac{2}{N})}{8C(n+4k)N^{n+4k-1}}(1 - \frac{\pi^{2}}{(\pi + \frac{\delta}{2})^{2}}), $$
and if
\[ \frac{1}{vol_{f}(M)}\int_{M} ((n-1) - Ric_{f-})_{+}e^{-f}dv_{g}< \epsilon, \]
then,
\begin{equation}
\sum_{i=1}^{n-1} I(bE_{i},bE_{i}) < - \frac{(n-1)L}{2}\left(1 - \frac{\pi^{2}}{L^{2}}\right) + \frac{(n-1)L}{2}\left(1 - \frac{\pi^{2}}{(\pi + \frac{\delta}{2})^{2}}\right). \label{eqn8}
\end{equation}
Because $\gamma $ is a minimal geodesic such that  
$$\sum_{i=1}^{n-1} I(bE_i, bE_i) \geq 0,$$
Then by (\ref{eqn8}), we obtain
$$
\frac{(n-1)L}{2}\left(1 - \frac{\pi^{2}}{(\pi + \frac{\delta}{2})^{2}}\right) \geq \frac{(n-1)L}{2}\left(1 - \frac{\pi^{2}}{L^{2}}\right). $$
This inequality gives
\begin{equation}
L \leq \pi + \frac{\delta}{2}. \label{eqn10}
\end{equation}
Finally, by the inequality (\ref{eqn5}), (\ref{eqn6}), and (\ref{eqn10}), we have
$$
{\rm diam}(M) = D < \pi + \delta.
$$
This completes the proof.
\section{Proof of Theorem~\ref{main2} and \ref{main3}}
We will prove Theorem~\ref{main2} in this section following the same setting for Theorem~\ref{main1}. Let $\gamma$ be a minimizing unit speed geodesic from $x \in \overline{A_{1}}$ to $y \in \overline{A_{2}}$ of length $L = d(x,y)$. Likewise, let $\{E_{1},\cdots, E_{n} = \gamma'\}$ be a parallel orthonormal frame along $\gamma$ and a smooth function $b \in C^{\infty}([0,L])$ such that $b(0) = b(L) = 0$.
For the proof, we need the following result. 

\begin{lemma} \label{lem1}
 Let $(M,g,e^{-f}dv_{g})$ be a complete Riemannian manifold with $Ric_{f} \geq (n-1)H$, $H < 0$, and $|f| \leq k$. Then, for any fixed $R > \pi$, there exists $\epsilon = \epsilon(n+4k,H,R,\delta)$ such that if
$$
\frac{1}{vol_{f}(B(p,R))} \int_{B(p,R)} ((n-1) - Ric_{f-})_{+} e^{-f}dv_{g} < \epsilon(n+4k,H,R,\delta)
$$
for some $B(p,R) \subset M$, then $M = B(p,R) \subset B(p,\pi + \delta)$.
\end{lemma}
\begin{pol}
 Let $sn_{H}(r)$ be the solution to 
$$sn_{H}'' + H sn_{H} = 0$$ 
such that $sn_{H}(0) = 0$ and $sn_{H}'(0) = 1,$ then 
$$m_{H}^{n} = (n-1)\frac{sn_{H}'}{sn_{H}}.$$ When $H < 0$, this solution is given by $sn_{H}(r) = \frac{1}{\sqrt{-H}}\sinh(\sqrt{-H}r)$.
By the inequality (\ref{eqn2}),
$$
m_{f}(r) \leq m_{H}(r) + \frac{4k sn_{H}'(r)}{sn_{H}(r)} = (n + 4k -1) \frac{sn_{H}'(r)}{sn_{H}(r)} = m_{H}^{n+4k}(r), 
$$
thus,
\be \frac{vol_{f}(B(p,R))}{vol_{f}(B(p,r))} \leq \frac{vol_{H}^{n+4k}(R)}{vol_{H}^{n+4k}(r)}. \label{eqn11}\ee
 
 If we set $p \in M$ and $W = B(p,R)$, then $q$ will be any point in $W$ satisfying $\pi + 4r < d(p,q) < R -3r$, where $0 < r < \frac{1}{8}(R - \pi)$ is to be determined, and $A_{1} = B(p,r), A_{2} = B(q,r)$. Thus, by (\ref{eqn11})
\bea 
\lefteqn{\int_{\gamma_{x,y}} ((n-1) - Ric_{f-})_{+} ds}\\
 &\leq & 2r C(n+4k,H,R)\left(\frac{1}{vol_{f}(B(p,r))} + \frac{1}{vol_{f}(B(q,r))}\right) \\
 & & \times \int_{B(p,R)} ((n-1) - Ric_{f-})_{+} e^{-f}dv_{g}, \\
& \leq & 2r C(n+4k,H,R)(\frac{vol_{H}^{n+4k}(R)}{vol_{H}^{n+4k}(r)}\frac{1}{vol_{f}(B(p,R))} \\
& & + \frac{vol_{H}^{n+4k}(2R)}{vol_{H}^{n+4k}(r)}\frac{1}{vol_{f}(B(q,2R))})\int_{B(p,R)} ((n-1) - Ric_{f-})_{+} dv_{g}, \\
&\leq & 2rC(n+4k,H,R)\left(\frac{vol_{H}^{n+4k}(R) + vol_{H}^{n+4k}(2R)}{vol_{H}^{n+4k}(r)}\right) \\
& & \times \frac{1}{vol_{f}(B(p,R))}\int_{B(p,R)} ((n-1) - Ric_{f-})_{+} e^{-f}dv_{g}.\eea

   Let $r = \frac{1}{4}\delta$, where $\delta < \frac{1}{2}(R - \pi)$. Then
\bea
\int_{\gamma_{x,y}} ((n-1) - Ric_{f-})_{+} ds &\leq &\frac{1}{2}\delta C(n+4k,H,R)\left(\frac{vol_{H}^{n+4k}(R) + vol_{H}^{n+4k}(2R)}{vol_{H}^{n+4k}(\frac{1}{4}\delta)}\right) \\
& & \times \frac{1}{vol_{f}(B(p,R))} \int_{B(p,R)} ((n-1) - Ric_{f-})_{+} e^{-f}dv_{g}. 
\eea

 By the second variation of $\gamma$, 
\bea
\sum^{n-1}_{i=1} I(bE_{i},bE_{i})& \leq & - \frac{(n-1)L}{2}(1 - \frac{\pi^{2}}{L^{2}}) + \int^{L}_{0} ((n-1) - Ric_{f-})_{+} ds, \\
& \leq & - \frac{(n-1)L}{2}(1 - \frac{\pi^{2}}{L^{2}}) + \frac{1}{2}\delta C(n+4k,H,R)(\frac{vol_{H}^{n+4k}(R)}{vol_{H}^{n+4k}(\frac{1}{4}\delta)} \\
& &+ \frac{vol_{H}^{n+4k}(2R)}{vol_{H}^{n+4k}(\frac{1}{4}\delta)}) \frac{1}{vol_{f}(B(p,R))}\int_{B(p,R)}((n-1) - Ric_{f-})_{+} e^{-f}dv_{g}. \eea

Setting
$$
\epsilon = \frac{(n-1)(\pi + \frac{1}{2}\delta)}{\delta C(n+4k,H,R)
}\left(1 - \frac{\pi^{2}}{(\pi + \frac{1}{2}\delta)^{2}}\right)\frac{vol_{H}^{n+4k}(\frac{1}{4}\delta)}{vol_{H}^{n+4k}(R) + vol_{H}^{n+4k}(2R)},
$$
we obtain
$$
\sum^{n-1}_{i=1} I(bE_{i},bE_{i}) < -\frac{(n-1)L}{2}\left(1 - \frac{\pi^{2}}{L^{2}}\right) + \frac{(n-1)(\pi + \frac{1}{2}\delta)}{2}\left(1 - \frac{\pi^{2}}{(\pi + \frac{1}{2}\delta)^{2}}\right).$$
Moreover, by the minimality of $\gamma$, we have
$$
\frac{(n-1)(\pi + \frac{1}{2}\delta)}{2}\left(1 - \frac{\pi^{2}}{(\pi + \frac{1}{2}\delta)^{2}}\right) \geq \frac{(n-1)L}{2}\left(1 - \frac{\pi^{2}}{L^{2}}\right), $$
implying that
$$
L \leq \pi + \frac{1}{2}\delta. $$
By the triangle inequality,
\be
 d(p,q) \leq \pi + \delta. \label{eqnw1}
\ee
We assumed that $\pi + 4r < d(p,q) < R - 3r$, or $\pi + \delta < d(p,q) < R - 3r$.
However, by (\ref{eqnw1}), no geodesic starting from $p$ of a length greater than $\pi + \delta$ can be length minimizing, which implies that $B(p,R) \subset B(p,\pi + \delta)$. If $R $ goes to infinity,  then $B(p,R)$ tends to $M$. Hence, we may conclude that $M = B(p,R) \subset B(p,\pi + \delta)$.
\end{pol}

 Now we can prove Theorem~\ref{main2}.\\ \\
\begin{pot1}
 Note that Lemma~\ref{lem1} shows Theorem~\ref{main2} for $R > \pi$. Thus, it suffices to prove the case when $R \leq \pi$.

   Let $R^{\prime} > \pi$ be fixed. Then for any $R \leq \pi$, there exists $N = N(H, R, R^{\prime})>0$, such that any $R^{\prime}$-ball in $M$ can be covered by $N$ or fewer $R$-balls, $B(x_{i},R)$, $1\leq i\leq N$. Subsequently,
\bea
\lefteqn{\frac{1}{vol_{f}(B(z,R'))} \int_{B(z,R')} ((n-1) - Ric_{f})_{+} e^{-f}dv_{g}}\nonumber\\
& \leq & N(H,R,R') \frac{1}{vol_{f}(B(z,R'))} \sup_{x_{i}} \int_{B(x_{i},R)} ((n-1) - Ric_{f-})_{+} e^{-f}dv_{g}, \nonumber\\
& \leq& N(H,R,R')\, \frac{vol_{H}^{n+4k}(R + R')}{vol_{H}^{n+4k}(R')}\, \frac{1}{vol_{f}(B(z,R + R'))} \nonumber \\
& &\times \sup_{x_{i}} \int_{B(x_{i},R)} ((n-1) - Ric_{f-})_{+} e^{-f}dv_{g}, \nonumber\\
&\leq & N(H,R,R')(\frac{vol_{H}^{n+4k}(R + R')}{vol_{H}^{n+4k}(R')}) \sup_{x_{i}}\frac{1}{vol_{f}(B(x_{i},R))}\nonumber \\
& & 
\times \int_{B(x_{i},R)} ((n-1) - Ric_{f-})_{+} e^{-f}dv_g. 
\eea
Hence, we can conclude that
\bea
\lefteqn{\sup_{x} \frac{1}{vol_{f}(B(x,R'))} \int_{B(x,R')}((n-1) - Ric_{f-})_{+}e^{-f}dv_g},\nonumber\\
&\leq& N(H,R,R')\, \frac{vol_{H}^{n+4k}(R + R')}{vol_{H}^{n+4k}(R')}\nonumber\\
 & &\times \, \sup_{x} \frac{1}{vol_{f}(B(x,R))}\int_{B(x,R)}((n-1) - Ric_{f-})_{+} e^{-f}dv_g. 
\eea
\end{pot1}

 Finally, let us prove Theorem~\ref{main3}.\\ \\
\begin{pot2}
   Let $(\tilde{M},\tilde{g})$ be a Riemannian universal cover of $(M,g)$. Because the inequality
$$\sup_{x} \frac{1}{vol_{f}(B(x,R))} \int_{B(x,R)} ((n-1) - Ric_{f-})_{+} e^{-f}dv_{g} < \tilde{\epsilon}(n+4k,H,R)$$
holds on $(M,g)$, the same inequality holds on $(\tilde{M},\tilde{g})$. 

   Based on this, it is easy to see that Theorem~\ref{main2} with $\delta =0$ also holds. When $R\leq \pi$, we just need to follow the proof of Theorem~\ref{main2}; moreover, when $R>\pi$, setting
$\pi <d(p,q) <R$ with $r=\frac {\pi}4$, we can prove that ${\rm diam}\, (M)\leq \pi$.
   Hence, we can conclude that $\tilde{M}$ is compact, implying that the fundamental group $\pi_{1}(M)$ is finite.
\end{pot2}

\end{document}